\documentclass{amsart}[11pt]

\usepackage{amsmath}
\usepackage{amsthm, amsfonts, latexsym, amssymb, mathrsfs}

\newcommand{\bC}{\mathbb{C}}
\newcommand{\C}{\bC}
\newcommand{\bF}{\mathbb{F}}

\newcommand{\Fq}{{\bF}_q}
\newcommand{\Fqbar}{\overline{\Fq}}

\newcommand{\bQ}{\mathbb{Q}}
\newcommand{\Q}{\bQ}
\newcommand{\Qbar}{\overline{\Q}}

\newcommand{\Ql}{{\bQ}_{\ell}}
\newcommand{\Qlbar}{\overline{\Ql}}

\newcommand{\bZ}{\mathbb{Z}}
\newcommand{\Z}{\bZ}

\newcommand{\Zl}{{\bZ}_{\ell}}

\newcommand{\calO}{\mathcal{O}}

\newcommand{\OP}[1]{\operatorname{#1}}
\newcommand{\Gal}{\OP{Gal}}
\newcommand{\Spec}{\OP{Spec}}

\newcommand{\cris}{\mathrm{cris}}
\newcommand{\et}{{\acute{\mathrm{e}}\mathrm{t}}}
\newcommand{\Frob}{\mathrm{Frob}}
\newcommand{\ord}{\mathrm{ord}}

\newcommand{\Map}{\longrightarrow}

{\theoremstyle{plain}\newtheorem{theorem1}[subsection]{Theorem}}
{\theoremstyle{plain}\newtheorem{theorem}[equation]{Theorem}}
{\theoremstyle{plain}\newtheorem{proposition}[equation]{Proposition}}
{\theoremstyle{plain}\newtheorem{lemma}[equation]{Lemma}}
{\theoremstyle{plain}}
{\theoremstyle{plain}\newtheorem{corollary}[equation]{Corollary}}

{\theoremstyle{definition}}

{\theoremstyle{remark}}
{\theoremstyle{remark}\newtheorem{remark1}[subsection]{Remark}}

\begin{document}

\title
[Symmetry and parity]
{Symmetry and parity in Frobenius action on cohomology}
\author{Junecue Suh}


\begin{abstract}
We prove that the Newton polygons of Frobenius on the crystalline cohomology of proper smooth varieties satisfy a symmetry that results, in the case of projective smooth varieties, from Poincar\'e duality and the hard Lefschetz theorem. As a corollary, we deduce that the Betti numbers in odd degrees of any proper smooth variety over a field are even (a consequence of Hodge symmetry in characteristic zero), answering an old question of Serre. Then we give a generalization and a refinement for arbitrary varieties over finite fields, in response to later questions of Serre and of Katz.
\end{abstract}

\maketitle

Let $X$ be a projective smooth variety over $\C$. Then the complex manifold associated with $X$ is a K\"ahler manifold, and its cohomology $H^r(X, \Q)$ is equipped with a pure Hodge structure:
$$
H^{r}(X, \Q) \otimes_{\Q} \C = \bigoplus_{p+q=r} H^{pq},
$$
where $H^{pq} = H^q(X, \Omega^p_X)$ satisfies $\overline{H^{pq}} = H^{qp}$. In particular, one has the Hodge symmetry
$$
h^{pq} = h^{qp}, \mbox{ where } h^{pq} := \dim_{\C} H^{pq},
$$
which implies that $b_r = \dim_{\Q} H^r(X, \Q)$ is even when $r$ is odd. This imposes a nontrivial condition on the topology of projective smooth varieties (or K\"ahler manifolds). For instance, it keeps the Hopf manifold $(\C^{2} - \{ 0 \})/ \Z$ (which is a compact complex manifold) from being a projective smooth variety, because it has $b_1 = 1$.

In \cite[\S 5]{deligneLefschetz}, Deligne constructs a pure Hodge structure on the cohomology of any \textit{proper} smooth variety over $\C$. Thus the Hodge symmetry and the evenness of odd-degree Betti numbers extend to the proper smooth case.

Now let $X$ be a projective smooth variety over a finite field $k$. For either the $\ell$-adic cohomology $H^{\ast}_{\et}(X\otimes_k \bar{k}, \Ql)$ for $\ell$ invertible in $k$ or the crystalline cohomology $H^{\ast}_{\cris}(X/W)_K$ (where $K$ denotes the fraction field of $W=W(k)$), Poincar\'e duality (\cite{deligneArcata} and \cite{berthelot1}) and the hard Lefschetz theorem (\cite{deligneWeil2} and \cite{katzmessing}) endow $H^r(X)$ with a perfect pairing that is symmetric when $r$ is even and alternating when $r$ is odd. In particular, the odd-degree cohomology groups are still even-dimensional. In the case of crystalline cohomology, we also get a symmetry in the Newton polygon (for every $r$).

In this article, we first show that these symmetry and parity statements extend to the proper smooth case (\textit{without} a perfect pairing), answering an old question of Serre. In response to later questions of Serre and of Katz, we then prove more general and refined statements concerning arbitrary varieties over finite fields. These may be considered as more concrete, observable consequences of conjectural properties of odd-weight motives in characteristic $p$.

\subsection*{Acknowledgments}
I express my gratitude to Luc Illusie, from whom I learned of the question on the parity of Betti numbers. I am also grateful to Jean-Pierre Serre, who suggested a generalization and a refinement (together with evidence for the latter, based on the resolution of singularities) of his original question. I sincerely thank Nicholas Katz for many helpful discussions, encouragement, and his question on the determinant of Frobenius. I also thank Shenghao Sun for his helpful comments.

\section{Theorems of Katz-Messing and of Gabber}
Throughout this section, let $k$ be a finite field with $q=p^e$ elements. Recall that an algebraic integer $\alpha$ is called a \textbf{$q^r$-Weil integer}, or a $q$-Weil integer of weight $r$, if for any embedding $\sigma: \Q(\alpha) \Map \C$, we have
\begin{equation}
\label{Weil_int}
\sigma(\alpha) \cdot \overline{\sigma(\alpha)} = q^r.
\end{equation}
Any Weil integer generates over $\Q$ either a totally real field or a CM field.
\begin{theorem1}[(Deligne)]
\label{thm_deligne}
Let $X$ be a proper smooth variety over $k$ and $\ell$ a prime number different from $p$. Then for any integer $r\ge 0$, all the eigenvalues of the (geometric) Frobenius acting on $H^r_{\et}(X\otimes_k \bar{k}, \Ql)$ are $q^r$-Weil integers, and the reversed characteristic polynomial
$$
P_{r, \et}(T) := \det \left( 1 - T\Frob_q: H^r_{\et}(X\otimes_k \bar{k}, \Ql) \right)
$$
has integer coefficients and is independent of the choice of $\ell \neq p$.
\end{theorem1}
This is proved in \cite{deligneWeil1} in the projective case, and follows from a much more general theorem in \cite{deligneWeil2} in the proper case. By using results from \cite{deligneWeil2} (including the pgcd theorem on Lefschetz pencils), Katz and Messing proved:
\begin{theorem1}[(Katz-Messing, {\cite{katzmessing}})]
\label{thm_km}
Let $X$ be a projective smooth variety over $k$. For any $r\ge 0$, the polynomial
$$
P_{r, \cris}(T) := \det \left( 1- T \Frob_q : H^r_{\cris}(X/W)_K \right),
$$
where $\Frob_q := F^e$ and $F$ is the crystalline Frobenius, has integer coefficients and is equal to $P_{r, \et}(T)$.
\end{theorem1}
Using similar ideas, Gabber proved:
\begin{theorem1}[(Gabber, {\cite{gabber}})]
\label{thm_gabber}
Let $X$ be a projective smooth variety over $k$. For all but finitely many primes $\ell\neq p$, the $\Zl$-cohomology group $H^{\ast}_{\et}(X \otimes_k \bar{k}, \Zl)$ is torsion-free.
\end{theorem1}
Gabber's theorem extends to the case of an arbitrary field $k$, by applying base change theorems in \'etale cohomology after finding a model over a finitely generated $\Z$-algebra.

These theorems extend to the proper smooth case. In \cite{chls}, Chiarellotto and Le Stum extend Theorem \ref{thm_km}, by first comparing the crystalline cohomology with the rigid cohomology, and then combining Poincar\'e duality with earlier results \cite{ch} on the weights appearing in the rigid cohomology of smooth varieties.

It is no doubt well known to the experts, but the extension of Theorem \ref{thm_gabber} does not seem to be recorded in the literature. We give a sketch of proof for both theorems, based on a simple application of de Jong's theory of alterations.
\begin{theorem1}
\label{thm_folklore}
The conclusions of both Theorems \ref{thm_km} and \ref{thm_gabber} hold true for any proper smooth variety $X$ over $k$.
\end{theorem1}
\begin{proof}
We may assume that $X$ is connected, hence integral. By applying Chow's lemma \cite[Lem 5.6.1]{ega2} and then \cite[Th. 4.1]{dejong}, we get a projective, surjective and generically finite morphism $\pi:Y\Map X$ from a projective smooth variety $Y$ over $k$; we may also assume that $\pi$ is generically \'etale.

By Poincar\'e duality, the pullback map $\pi^{\ast}: H^i_{\cris}(X/W)_K \Map H^i_{\cris}(Y/W)_K$ defines $\pi_{\ast}: H^i_{\cris}(Y/W)_K \Map H^i_{\cris}(X/W)_K$. Then we have $\pi_{\ast}\circ \pi^{\ast} = \deg(\pi)$. Indeed, by the projection formula and the compatibility of the cycle class map with the proper push-forward, we have
$$
\pi_{\ast} \pi^{\ast} x = \pi_{\ast} ( [Y] \cdot \pi^{\ast} x ) = \pi_{\ast} ( [Y] ) \cdot x = \deg(\pi) [X] \cdot x = \deg(\pi) x,
$$
for any $x\in H^i_{\cris}(X/W)_K$. See \cite{gm} for the definition of the cycle class map and the relevant facts.

This implies that $\pi^{\ast}$ is an injection, so all the eigenvalues of $\Frob_q$ on $H^i_{\cris}(X/W)_K$ are $q^i$-Weil integers. Using the cohomological interpretation of the zeta function:
\begin{equation}
\label{coh-exp}
Z(X/\Fq, T) = \prod_{i=0}^{2 \dim X} P_{i, \cris} (T)^{(-1)^{i+1}} ,
\end{equation}
we recover a given $P_{i, \cris}$ as the weight $i$ part (\textit{cf.} the footnote in \S \ref{sec-review}) of the zeta function itself, i.e., the eigenvalues of $\Frob_q$ on $H^i_{\cris}(X/W)_K$ for $i$ odd (resp. even) are precisely the reciprocal zeroes (resp. poles) of the zeta function that are $q^i$-Weil integers. Because the zeta function lies in $\Q(T)$
and the notion of a $q^i$-Weil integer is $\Gal(\Qbar/\Q)$-invariant, it follows that the $P_{i, \cris}$ have coefficients in $\Q$, hence in $\Z$ since Weil integers are algebraic integers. Either by directly invoking Deligne's Theorem \ref{thm_deligne} in the proper smooth case or by applying the same $\pi_{\ast} \circ \pi^{\ast} = \deg (\pi)$ argument in $\ell$-adic cohomology for any $\ell \neq p$, we get the same characterization of the $P_{i, \et}$, and we have $P_{i, \cris} = P_{i, \et}$.

We also use the map $\pi$ to extend Gabber's theorem. For any ring $A$, put $H(X, A) = H^{\ast}_{\et}(X\otimes_k \bar{k}, A)$ and likewise for $H(Y,A)$. It is enough to show that the Gysin map $\pi_{\ast} : H(Y, \Zl) \Map H(X, \Zl)$ satisfies $\pi_{\ast} \circ \pi^{\ast} = \deg(\pi)$ for all $\ell \neq p$. For then $H(X, \Zl)$ will be torsion-free whenever $\ell \nmid p\deg(\pi)$ and $H(Y, \Zl)$ is torsion-free. Using Poincar\'e duality with $\Z/\ell^n$-coefficients and the normalisation \cite[XVIII, Th 2.9 (Var4)]{sga4}, one defines a projective system of maps $\pi_{\ast}: H(Y, \Z/\ell^n) \Map H(X, \Z/\ell^n)$ satisfying $\pi_{\ast} \circ \pi^{\ast} = \deg (\pi)$ on $H(X, \Z/\ell^n)$. One concludes by taking the limit.

\end{proof}
\begin{remark1}
The same applications to cycles as in \cite[Th. 2]{katzmessing} work in the proper smooth case, see \cite[3.5(a)]{illusie}.
\end{remark1}

\numberwithin{equation}{subsection}

\section{Symmetry and parity for proper smooth varieties}
After making some definitions about multisets, we state and prove our main theorem for proper smooth varieties. It turns out that the symmetry in the Newton polygon follows from Theorem \ref{thm_folklore} and a simple Galois theory argument. To show that the Betti numbers in odd degrees are even, however, we \textit{do} need the underlying $F$-isocrystal.

\subsection{Multisets}
\label{sec-multiset}
By a \textbf{multiset} (with finite support) $S$ of elements in an ambient set $\Sigma$, we mean a finite subset of $\Sigma$ with multiplicities. Formally, these multisets correspond to functions $\mu_S : \Sigma \Map \Z_{\ge 0}$ such that $\mu_S(\sigma)= 0$ for all but finitely many $\sigma\in \Sigma$.

Some notions of usual subsets generalize naturally to multisets. We say that $\sigma\in S$ if $\mu_S(\sigma)>0$, and that two multisets $S_1$ and $S_2$ in $\Sigma$ are disjoint if we have $\mu_{S_1}(\sigma) \mu_{S_2}(\sigma) = 0$ for all $\sigma\in \Sigma$. The cardinality $|S|$ of $S$ is defined as $\sum_{\sigma\in \Sigma} \mu_S(\sigma)$. If $\Sigma$ is an abelian group (such as $\Q$) and $r \in \Sigma$, we denote by $r - S$ the multiset $\{ r - s : s\in S \}$, counted with multiplicities; formally, $\mu_{r - S}(\sigma) = \mu_S(r-\sigma)$ and we say that $S$ is \textbf{$r$-autodual} if $S = r - S$. Finally, if a group $G$ acts on the set $\Sigma$, we denote by $gS$ the multiset $\{ gs : s\in S\}$, and we say that $S$ is $G$-invariant if $gS = S$ for all $g\in G$.

Let $\Qbar \subseteq \C$ be the algebraic closure of $\Q$ in $\C$, and let $v: \Qbar^{\times} \Map \Q$ be a $p$-adic valuation. If $S$ is a multiset in $\Qbar^{\times}$, we denote by $v(S)$ the multiset in $\Q$ that is the ``image'' of $S$ under $v$, counted with multiplicities. If $S$ is invariant under $\Gal(\Qbar / \Q)$, then $v(S) = v'(S)$ for any $p$-adic valuation $v'$ such that $v(p) = v'(p)$.

\subsection{Main theorem for proper smooth varieties}
\label{sec-main-prsm}
\begin{proposition}
\label{prop-symmetry}
Let $q = p^{e}$ be a power of a prime $p$ with $e\ge 1$, $r$ a nonnegative integer, $v$ a $p$-adic valuation on $\Qbar^{\times}$ normalized by $v(q) = 1$, and $f(T) \in 1+ T\Q[T]$ a polynomial. Write
$$
f(T) = \prod_{\beta\in S} (1-\beta T) \mbox{ ( or, formally, } \prod_{\beta\in \Qbar^{\times}} (1-\beta T)^{\mu_S(\beta)} \mbox{ )},
$$
where $S$ is the $\Gal(\Qbar / \Q)$-invariant multiset in $\Qbar^{\times}$ consisting of the reciprocal roots of $f(T)$.

Assume that every $\beta \in S$ is a $q^r$-Weil integer. Then $v(S)$ is $r$-autodual, i.e., $v(S) = r- v(S)$.
\end{proposition}
\begin{proof}
By (\ref{Weil_int}), complex conjugation acts as $\beta \mapsto q^r/\beta$ on $S$.
\end{proof}

\begin{theorem}
\label{thm_main}
Let $X$ be a proper smooth variety over a finite field $k$ of characteristic $p$, and let $r\ge 0$ be an integer. Then the multiset of slopes of Frobenius on $H^r_{\cris}(X/W)_K$ is $r$-autodual, i.e., for any $s \in [0, r]$, $s$ and $r-s$ appear with the same multiplicity. If $r$ is odd, then $\dim_K H^r_{\cris}(X/W)_K$ ($=\dim_{\Ql} H^r_{\et}(X\otimes_k \bar{k}, \Ql)$ for any $\ell \neq p$) is even.
\end{theorem}
\begin{proof}
The first statement follows from applying Proposition \ref{prop-symmetry} to the polynomial $P_{i, \cris}(T)$; Theorem \ref{thm_folklore} verifies the assumption made in the proposition. For the second statement, note that the multiplicity of $r/2$ as a slope of the $F$-isocrystal $H^r_{\cris}(X/W)_K$ is necessarily even, from Dieudonn\'e-Manin classification. (The key point here is that the ``$\ord_p$-slopes'' of the $\sigma$-linear $F$ on an $F$-isocrystal over $K=W(\Fq)[1/p]$ are equal to the ``$\ord_q$-slopes'' of the eigenvalues of the $K$-linear $\Frob_q$, and that in any $F$-isocrystal, the multiplicity of a slope, written in lowest terms, is always an integral multiple of its denominator.)
\end{proof}
We remark that the use of the underlying $F$-crystal is essential in our proof: If $q$ is an even power of a prime and $r\ge 1$ is any odd integer, the linear polynomial $f(T) = 1 - \sqrt{q^r} T$ satisfies the symmetry, but not the parity.

By a standard ``spreading out'' argument, the parity statement extends to:
\begin{corollary}
\label{cor-parity}
Let $X$ be a proper smooth variety over any field $k$, and let $K$ be a separably closed extension of $k$. Then for any odd integer $r\ge 1$, the $r$-th Betti number
$$
b_{r, \ell} := \dim_{\Ql} H^r (X\otimes_k K, \Ql)
$$
is independent of $\ell \neq \mathrm{char}(k)$ and is even.
\end{corollary}
This answers the parity question that we learned of from Illusie. He says it was originally Serre who asked him the question. It also appeared in print on p. 394 of \cite{deligneWeil2}.

The theorem also extends to crystalline cohomology over more general fields (\textit{cf.} \cite{katz}):
\begin{corollary}
Theorem \ref{thm_main} is valid for the $F$-isocrystal underlying the crystalline cohomology of a proper smooth variety over any perfect field of characteristic $p>0$.
\end{corollary}

Even in this proper and smooth case, one interesting question remains unanswered by Theorem \ref{thm_main}. When $X$ is projective and smooth and $r$ is odd, the presence of a Frobenius-equivariant, alternating and nondegenerate pairing on $H^r(X)$ forces the determinant of Frobenius on it to be $q^{r(\dim H^r(X))/2}$. In the proper smooth case, Theorem \ref{thm_main} shows this \textit{only up to sign}. This question of sign will be settled later in \S \ref{sec-signs}, see Corollary \ref{question_katz}.

\section{Generalization and refinement} \label{sec-gen-question}
When he received an earlier draft containing Theorem \ref{thm_main}, Serre raised a generalized version of his original question, concerning arbitrary varieties over finite fields, as well as a refinement. We give answers below, in Theorems \ref{thm_gen} and \ref{thm_refine}.

Throughout this section, $X$ will be a separated scheme of finite type over a finite field $k=\Fq$ of characteristic $p>0$. We denote by $\Qbar$ the algebraic closure of $\Q$ in $\C$, and $\Q^{\mathrm{cm}}$ the compositum of all CM fields in $\Qbar$. We fix an embedding of $\Qbar$ into a chosen algebraic closure of $K=W(k)[1/p]$ and get a $p$-adic valuation $v$ on $\Qbar$, normalized by $v(q) = 1$. The \textbf{slope} of $\alpha \in {\Qbar}^{\times}$ will mean $v(\alpha)$. For $\ell\neq p$, we regard $\Qbar$ as a subfield of $\Qlbar$ via a fixed embedding.

\subsection{Review}
\label{sec-review}
For any prime $\ell\neq p$, we write $H^i_{\ell}= H^i_{\et, c}\left( X\otimes_k \bar{k}, \Ql \right)$ for the $\ell$-adic cohomology with compact support. By \cite[\S 3.3]{deligneWeil2}, every eigenvalue of Frobenius on $H^i_{\ell}$ is a $q$-Weil integer of some weight (see \cite[4.3]{illusie} for integrality). One expects that the polynomial
$$
P_{i, \ell}(T) = \det \left( 1-T\Frob_q : H^i_{\ell} \right)
$$
lies in $1+T\Z[T]$ and is independent of $\ell$, but neither the integrality nor the independence is known in general.

Write $H^i_p = H^i_{c, \mathrm{rig}}(X/K)$, the rigid cohomology with compact support \cite{berthelot2} of $X/k$. By cohomological descent in rigid cohomology \cite{cht}, \cite{tsuzuki} applied to proper hypercoverings obtained from alterations \cite{dejong}, one also proves that every eigenvalue of Frobenius on $H^i_p$ is a $q$-Weil integer; \textit{cf.} \cite[\S 5.2]{tsuzuki}. Again, one expects, but does not know in general, that the corresponding polynomial $P_{i,p}(T)$ lies in $1+T\Z[T]$ and coincides with $P_{i, \ell}(T)$ for any $\ell\neq p$.

For \textit{every} prime $\ell$ (including $\ell = p$), put $P_{\ast, \ell}(T):= \prod_i P_{i, \ell}(T)$, the reversed characteristic polynomial of $\Frob_q$ on $H^{\ast}_{\ell} = \oplus_i H^{i}_{\ell}$, and write
\begin{equation}
\label{theform}
P_{\ast, \ell}(T) = \prod_{\beta\in S^{\ast}_{\ell}} (1-\beta T),
\end{equation}
where $S^{\ast}_{\ell}$ is a multiset in $(\Q^{\mathrm{cm}})^{\times}$.

For each integer $r\ge 0$ and each prime $\ell$, the sum, over an algebraic closure of $\Ql$ or $K$, of the generalized eigenspaces in $H^i_{\ell}$ (resp. $H^{\ast}_{\ell}$) of $\Frob_q$ with $q^r$-Weil integer eigenvalues descends to a subspace $H^{i, (r)}_{\ell}$ of $H^{i}_{\ell}$ (resp. $H^{\ast, (r)}_{\ell}$ of $H^{\ast}_{\ell}$). They are the kernels of polynomials in $\Frob_q$ with $\Z$-coefficients (depending on $r$ and $\ell$), and $H^{i, (r)}_p$ and $H^{\ast, (r)}_p$ are sub-$F$-isocrystals. Let $S^{i, (r)}_{\ell}$ and $P^{(r)}_{i, \ell}$ (resp. $S^{\ast, (r)}_{\ell}$ and $P^{(r)}_{\ast, \ell}(T)$) be the corresponding (as in (\ref{theform})) multiset and polynomial, respectively. The previous expectations lead us to expect that $S^{i, (r)}_{\ell}$ and hence $S^{\ast, (r)}_{\ell}$ should be independent of $\ell$.

Denote by $\mu_{\ell}^{(r)}$ the multiplicity function (see \S \ref{sec-multiset}) of $S^{\ast, (r)}_{\ell}$. We do not know that $\mu_{\ell}^{(r)}$ is independent of $\ell$, nor that $\mu_{\ell}^{(r)}$ takes a constant value on any $\Gal(\Q^{\mathrm{cm}}/ \Q)$-orbit, but we \textit{do know} that $\mu_{\ell}^{(r)} \mod 2$ satisfies these properties. This follows from the cohomological interpretation of the zeta function:
$$
\prod_i P_{i, \ell}(T)^{(-1)^{i+1}} = Z(X/\Fq, T)
$$
(see \cite{deligneArcata} for $\ell \neq p$ and \cite{els} for $\ell = p$), which implies
\begin{equation}
\prod_i P^{(r)}_{i, \ell}(T)^{(-1)^{i+1}} =  Z^{(r)}(X/\Fq, T),
\label{coh-int-zeta1}
\end{equation}
where the right hand side is the weight $r$ part\footnote{
Write
$$
Z(X/\Fq, T) = \prod_{\alpha \in A} (1- \alpha T) / \prod_{\beta \in B} (1- \beta T)
$$
with disjoint multisets $A$ and $B$ in $(\Q^{\mathrm{cm}})^{\times}$, collect the $q^r$-Weil integers in $A$ and $B$ into $A(r)$ and $B(r)$, and define
\begin{equation}
\label{coh-int-zeta2}
Z^{(r)} (X/\Fq, T) = \prod_{\alpha \in A(r)} (1-\alpha T) / \prod_{\beta \in B(r)} (1-\beta T) \in \Q(T).
\end{equation}
}
of $Z(X/\Fq, T)$. In particular, the parity of the cardinality:
$$
\left| S^{\ast, (r)}_{\ell} \right| = \sum_{\beta\in (\Q^{\mathrm{cm}})^{\times}} \mu_{\ell}^{(r)} (\beta)
$$
is independent of $\ell$.
\subsection{Parity and symmetry for general varieties}
\begin{theorem}
\label{thm_gen}
Let $X$ be a separated scheme of finite type over $\Fq$ and let $r\ge 1$ be an odd integer.
Then the degree (as a rational function) of the weight $r$ part of the zeta function $Z^{(r)}(X/\Fq, T)$ is even, and for any $\ell$, the cardinality $|S^{\ast, (r)}_{\ell}|$ ($=$ the number of $q^r$-Weil integers, counted with multiplicities, occuring as Frobenius eigenvalues in the total cohomology ($\ell$-adic or rigid, according as $\ell\neq p$ or $\ell=p$) with compact support) is also even.
\end{theorem}
\begin{proof}
Write (\ref{coh-int-zeta2}) in the reduced form:
\begin{equation}
Z^{(r)}(X/\Fq, T) = \frac{f(T)}{g(T)} (1 - \sqrt{q^r} T)^{m_0} (1 + \sqrt{q^r} T)^{m_1},
\label{reducedform}
\end{equation}
where $f(T)$ and $g(T)$ are relatively prime polynomials in $1+T\Z[T]$, of which neither of $\pm \sqrt{q^r}$ is a reciprocal root, and $m_0, m_1 \in \Z$. Here we isolate the cases of $\pm\sqrt{q^r}$, because these are the fixed points of complex conjugation acting on the set of $q^r$-Weil integers.

The degrees of both $f(T)$ and $g(T)$ are even, because complex conjugation acts without fixed points on their sets of reciprocal roots, so it remains to prove that $m_0+m_1$ is even. To see this, note that the multiplicity of $r/2$ in the slopes of the reciprocal roots of $f(T)$ (resp. of $g(T)$) is necessarily even, again because complex conjugation acts as $\beta \mapsto q^r/\beta$ without fixed points on the reciprocal roots. By the classification of Dieudonn\'e-Manin, the multiplicity of $r/2$ as slope in the $F$-isocrystal $H^{\ast, (r)}_p$ is even. These two facts imply, by (\ref{coh-int-zeta1}), that $m_0+m_1$ is even.
\end{proof}
One can also show that the multiset of Frobenius slopes in $H^{\ast, (r)}_p$ (for any $r\ge 0$) is $r$-autodual modulo $2$, i.e., the multiplicity of a slope $s$ and that of $r-s$ in $H^{\ast, (r)}_p$ are either both even or both odd.

When $X/\Fq$ is smooth, the statement of Theorem \ref{thm_gen} holds with ordinary cohomology in place of cohomology with compact support, by Poincar\'e duality.

\subsection{Signs}
\label{sec-signs}
\begin{theorem}
\label{thm_refine}
Let $X/\Fq$ be a separated scheme of finite type and let $r\ge 1$ be an odd integer. Then the multiplicity of $\sqrt{q^r}$ as a reciprocal root or pole in $Z(X/\Fq, T)$ is even (hence so is the multiplicity of $\sqrt{q^r}$ as Frobenius eigenvalue in the total cohomology with compact support, as in Theorem \ref{thm_gen}). The same is true for $-\sqrt{q^r}$.
\end{theorem}
The second statement follows from the first, in view of the proof of Theorem \ref{thm_gen}.

For $X/\Fq$ and $r$ as above, define $m(X/\Fq, r) \in \Z$ as the order of zero or pole at $T=1/\sqrt{q^r}$ of the zeta function $Z(X/\Fq, T)$, or equivalently that of $Z^{(r)}(X/\Fq, T)$. So we need to prove that $m(X/\Fq, r)$ is even.

First we give some preliminary lemmas.

\begin{lemma}
Suppose that $X$ is a separated scheme of finite type over a finite extension $\bF_{q^e}$ of $\Fq$, and let $X_0/\Fq$ be $X$ viewed as an $\Fq$-scheme. Then we have
$$
m(X/\bF_{q^e}, r) = m(X_0/\Fq, r).
$$
\label{respect_constants}
\end{lemma}
\begin{proof}
By definition of the zeta function of varieties over finite fields, we have
$$
Z(X_0/\Fq, T) = Z(X/\bF_{q^e}, T^e).
$$
Write $Z(X/\bF_{q^e}, T) = f(T)\cdot (1-\sqrt{q^{er}} T)^m$, where $m\in \Z$ and $f(T)\in \Q(\sqrt{q^{er}})(T)$ is defined and takes a nonzero value at $T=1/\sqrt{q^{er}}$. By definition, $m=m(X/\bF_{q^e}, r)$. On the other hand, we get
$$
Z(X_0/\Fq, T) = f(T^e) \cdot \left( 1-(\sqrt{q^r} T)^e \right)^m.
$$
Since $f(T^e)$ is defined and nonzero at $T=1/\sqrt{q^r}$, it follows from the cyclotomic factorization of the last factor that $m=m(X_0/\Fq, r)$.
\end{proof}

\begin{lemma}
Let $U\subseteq X$ be an open subset with complement $F$. Then we have
$$
m(X/\Fq, r) = m(U/\Fq, r) + m(F/\Fq, r).
$$
In particular, if two of the three are even, then so is the third.
\label{choppingup}
\end{lemma}
\begin{proof}
It follows from the definition of the zeta function that
$$
Z(X/\Fq, T) = Z(U/\Fq, T) Z(F/\Fq, T).
$$
\end{proof}

\begin{lemma}
Suppose that $X/\Fq$ is connected, projective and smooth and that $G$ is a finite group of automorphisms acting $\Fq$-linearly on $X$. Denote by $Y=X/G$ the quotient scheme. Then for any odd integer $r\ge 1$, both $m(X/\Fq, r)$ and $m(Y/\Fq, r)$ are even.
\label{rationalhomology}
\end{lemma}

\begin{proof}
Choose an auxiliary prime $\ell\neq p$, and an ample line bundle $L$ on $X$. By replacing $L$ with its $G$-norm ($=\otimes_{g\in G} g^{\ast}L$) if necessary, we may assume that the cohomology class of $L$ is fixed by $G$. By the hard Lefschetz theorem and Poincar\'e duality, $H^r(X) = H^r(X \otimes_{\Fq} \Fqbar, \Ql)$ has an alternating nondegenerate pairing $\langle \cdot , \cdot \rangle_L: H^r(X) \times H^r(X) \Map \Ql(-r)$ that is $G$-invariant and Frobenius-equivariant. Moreover, we have
$$
H^i(Y\otimes_{\Fq} \Fqbar, \Ql) = H^i(X\otimes_{\Fq} \Fqbar, \Ql)^G \mbox{ for any } i\ge 0.
$$
The $G$-invariance of $\langle \cdot , \cdot \rangle_L$ implies that it restricts to a nondegenerate alternating pairing on $H^r(Y)$, and Frobenius-equivariance implies that the multiplicity of $\sqrt{q^r}$ (in each of $H^r(X)$ and $H^r(Y)$) as an eigenvalue of Frobenius is even.
\end{proof}

Now let $X/\Fq$ be separated of finite type and $r\ge 1$ an odd integer. By a repeated use of Lemma \ref{choppingup}, we may assume that $X$ is integral, normal and projective over $\Fq$.

The ring $\Gamma:=\Gamma(X, \calO_X)$ is a finite field extension of $\Fq$. By Lemma \ref{respect_constants}, we may replace $\Fq$ by $\Gamma$ and assume that $X/\Fq$ is geometrically connected (hence geometrically integral, given that normality over a perfect field implies geometric normality \cite[Prop. 6.7.4]{ega4}).

We proceed by induction on dimension, and assume that we know $m(Z/\Fq, r)$ is even for every $Z$ of dimension strictly less than $\dim X$. By the induction hypothesis and Lemma \ref{choppingup}, it suffices to find a nonempty open $U \subseteq X$ for which $m(U/\Fq, r)$ is even.

By \cite[Th. 5.13]{dejong2} (applied to $S=\Spec(\Fq)$ with $G=\{ 1\}$ in the notation of \textit{loc. cit.}), there exist (A) a connected, projective and smooth scheme $X'$ over $\Fq$ with an action of a finite group $G$ and (B) a proper, surjective, generically finite and $G$-invariant morphism $\pi: X'\Map X$ such that the field extension $\Fq(X) \subseteq (\Fq(X'))^G$ is purely inseparable.

Let $X''= X'/G$ be the quotient scheme with the induced map $\pi'':X''\Map X$. There exists a nonempty open subset $U \subseteq X$ such that the restriction of $\pi''$ to $U'':= (\pi'')^{-1}(U)$ is finite and flat. By the condition on the function field extension in (B) above, we may assume that the restriction $\pi''|_{U''}: U''\Map U$ is a universal homeomorphism.
Since $U''$ is a dense open subset of $X''$, Lemmas \ref{rationalhomology} and \ref{choppingup} plus the induction hypothesis show that $m(U/\Fq, r)=m(U''/\Fq, r)$ is even. This completes the proof of Theorem \ref{thm_refine}.

\begin{corollary}
\label{question_katz}
Let $X$ be a proper smooth variety over $\Fq$, and let $r\ge 1$ be an odd integer. Then the determinant of $\Frob_q$ on the $r$-th cohomology (either $\ell$-adic or crystalline) is equal to $q^{r b_r/2}$, where $b_r=b_r(X/\Fq)$ is the dimension of the cohomology.
\end{corollary}
This answers the question at the end of \S \ref{sec-main-prsm}, raised independently by Serre and Katz.

We note that, as Serre kindly pointed out to us, a more straightforward proof of Theorem \ref{thm_refine} follows from the resolution of singularities in characteristic $p>0$ (which isn't available yet): We could express the zeta function of any variety in terms of those of projective smooth varieties directly, without the complications that had to be dealt with in our proof using equivariant alterations.

\end{document}